\documentclass[11pt]{article}
\usepackage[a4paper,margin=25mm]{geometry}
\usepackage[T1]{fontenc}
\usepackage{lmodern}
\usepackage{amsmath,amssymb,amsthm,mathtools,microtype}
\usepackage[hidelinks]{hyperref}
\hypersetup{pdftitle={Triangles in C5-free graphs}}
\newtheorem{theorem}{Theorem}[section]
\newtheorem{proposition}[theorem]{Proposition}
\newtheorem{lemma}[theorem]{Lemma}
\newcommand{\ex}{\operatorname{ex}}
\newcommand{\tr}{\operatorname{tr}}
\newcommand{\op}{\mathrm{op}}
\newcommand{\Frob}{\mathrm F}
\allowdisplaybreaks[1]
\title{Triangles in $C_5$-free graphs}
\author{Guorong Gao\thanks{School of Mathematics and Statistics,
    Fuzhou University, Fuzhou 350108, China.
    Email \href{mailto:grgao@fzu.edu.cn}{\nolinkurl{grgao@fzu.edu.cn}}.} \and
  Jialin He\thanks{School of Mathematical Sciences, Key Laboratory of MEA,
    and Shanghai Key Laboratory of PMMP, East China Normal University,
    Shanghai 200241, China.
    Email \href{mailto:jlhe@math.ecnu.edu.cn}{\nolinkurl{jlhe@math.ecnu.edu.cn}}.} \and
  Tianchi Yang\thanks{School of Mathematics, Georgia Institute of Technology,
    Atlanta, GA 30332, USA.
    Email \href{mailto:tyang439@gatech.edu}{\nolinkurl{tyang439@gatech.edu}}.}}
\date{}
\begin{document}
\maketitle
\vspace{-1em}
\begin{abstract}
We determine asymptotically the maximum number of triangles in a
$C_5$-free graph. More precisely, we prove that
\(
 \ex(n,K_3,C_5)=\left(\frac1{3\sqrt3}+o(1)\right)n^{3/2},
\)
establishing the asymptotic optimality of a construction of
Bollob\'as and Gy\H{o}ri. The proof combines counts of paths of length
two with spectral estimates for adjacency matrices and graph removal lemmas.
\end{abstract}

\section{Introduction}
For graphs $F$ and $H$, let $\ex(n,F,H)$ denote the maximum number
of unlabeled copies of $F$ in an $n$-vertex $H$-free graph. All graphs in this
paper are finite and simple, and copies are not required to be
induced. We study $\ex(n,K_3,C_5)$, the maximum number of triangles
in a graph containing no cycle of length five.

Bollob\'as and Gy\H{o}ri \cite{BG} initiated the study of this
problem, motivated by a question of Erd\H{o}s \cite{Erdos84} on the
maximum number of pentagons in an $n$-vertex triangle-free graph. Erd\H{o}s
conjectured that this number is at most $(n/5)^5$; the conjecture
was proved by Grzesik \cite{Grzesik} and, independently, by
Hatami, Hladk\'y, Kr\'a\v{l}, Norine and Razborov \cite{HHKNR}.
Reversing the roles of triangles and pentagons, Bollob\'as and
Gy\H{o}ri proved that
\[
 \left(\frac1{3\sqrt3}+o(1)\right)n^{3/2}
 \le \ex(n,K_3,C_5)
 \le \left(\frac54+o(1)\right)n^{3/2}.
\]
Their lower-bound construction starts with a bipartite $C_4$-free
graph, replaces each vertex in one part by two adjacent copies,
and retains the original adjacencies at both copies. Applying this
construction to the incidence graph of a finite projective plane
gives the constant $1/(3\sqrt3)$.

The upper bound has since been improved several times. Alon and
Shikhelman \cite{AS} reduced the leading constant to $\sqrt3/2$.
Ergemlidze, Methuku, Salia and Gy\H{o}ri \cite{EGMS} obtained
$1/(2\sqrt2)$ by describing the blocks formed by edge-sharing
triangles and selecting a $\{C_4,C_5\}$-free subgraph with one edge
for each triangle. They then applied a theorem of Erd\H{o}s and
Simonovits \cite{ES}. Ergemlidze and Methuku \cite{EM} improved
the constant to $1/(3\sqrt2)$ by counting paths of length five,
and further to approximately $0.231975$ through a connection with
hypergraphs of girth six. Their estimates, together with the sparse
triangle-removal theorem of Conlon, Fox, Sudakov and Zhao
\cite{CFSZ}, also yield the constant $2^{-9/4}$. Here the removal theorem shows that only $o(n^{3/2})$ triangles
belong to copies of $K_4$, so these triangles contribute only an
$o(n^{3/2})$ error in the estimates from Section~2 of \cite{EM}. Lv, He and Lu \cite{LHL} subsequently
proved the upper bound $(1/(2\sqrt6)+o(1))n^{3/2}$.
Our main result proves that the leading constant in the
Bollob\'as--Gy\H{o}ri construction is optimal.

\begin{theorem}\label{thm:main}
\(
 \ex(n,K_3,C_5)=\left(\frac1{3\sqrt3}+o(1)\right)n^{3/2},
\)
as $n\to\infty$.
\end{theorem}

This problem is an instance of the generalized Tur\'an problem,
developed systematically by Alon and Shikhelman \cite{AS}; see
Gerbner and Palmer \cite{GP} for a survey. For longer odd cycles,
Gy\H{o}ri and Li \cite{GL} related the number of triangles in a
$C_{2k+1}$-free graph to extremal numbers for even cycles.
F\"uredi and \"Ozkahya \cite{FO} refined these connections in their
work on hypergraph cycles, obtaining, in particular, another proof
that $\ex(n,K_3,C_5)=O(n^{3/2})$.

\paragraph{Overview of the proof.}
By the decomposition in \cite{EGMS} and the sparse
triangle-removal theorem of \cite{CFSZ}, we can discard $o(n^{3/2})$
triangles so that the remaining triangles form edge-disjoint books.
In a book, the triangles share a common edge, called the base.
Deleting the bases leaves the \emph{spoke graph}, whose edges are
called \emph{spokes}. Each surviving triangle contributes exactly
two spokes, so it suffices to bound their number.

The spoke graph of the Bollob\'as--Gy\H{o}ri construction is bipartite.
For any bipartite spoke graph, counting second neighbors on the two
sides gives the sharp upper bound. For the general case, suppose that the
triangle count exceeds the proposed bound by a fixed multiple of
$n^{3/2}$ and that every vertex lies in sufficiently many triangles.
We group the vertices into classes. After deleting class pairs that
account for few spokes, the graph on the classes is bipartite. Its
two sides give a partition of the original vertices with few spokes
inside the parts. Discarding the triangles that contain these spokes
leaves more triangles than the bipartite bound permits. To remove
the minimum triangle-degree assumption, we delete vertices of small
triangle degree and sum the losses.

\paragraph{Organization and notation.}
Section~\ref{sec:framework} gives the book reduction and states the
main lemmas. We prove Theorem~\ref{thm:main} in Section~\ref{sec:proof}.
Section~\ref{sec:bipartite} proves the bipartite bound
(Lemma~\ref{lem:bipartite}). Section~\ref{sec:matrix} constructs the
finite weighted approximation of Lemma~\ref{lem:model}, which we use
in Section~\ref{sec:reduction} to prove the approximate bipartition
in Lemma~\ref{lem:near}.
All asymptotic notation refers to $n\to\infty$, with the stated
parameters fixed. Unless otherwise specified, error terms are
uniform over the book families under consideration.

\section{Proof of Theorem~\ref{thm:main}}\label{sec:framework}
To prove the upper bound, we first pass from a $C_5$-free graph to
a family of books. We then state the lemmas about its spokes and
use them to prove Theorem~\ref{thm:main}.

\subsection{Reduction to book families}\label{sec:books}
Let $G$ be an $n$-vertex $C_5$-free graph. Delete the edges that lie
in no triangle. Every remaining edge now lies in a triangle, and
the triangle count is unchanged.

A \emph{book} is a graph consisting of triangles $stp_1,\ldots,stp_k$,
where $k\ge 1$ and $s,t,p_1,\ldots,p_k$ are distinct vertices.
We call $st$ its \emph{base edge}, or simply its \emph{base}, and
the vertices $p_1,\ldots,p_k$ its \emph{pages}. The edges $sp_i$ and
$tp_i$ are its \emph{spokes}.

Two triangles belong to the same \emph{triangle block} if they can
be joined by a sequence of triangles in which consecutive triangles
share an edge. We also use the term for the union of the edges in
such a block. The following proposition follows from Claims~1 and~2
in Section~2 of \cite{EGMS}.

\begin{proposition}[\cite{EGMS}]\label{prop:blocks}
Let $G$ be a $C_5$-free graph.
\begin{enumerate}
\item[(i)] Every triangle block of $G$ is a book or a $K_4$.
\item[(ii)] If every edge of $G$ lies in a triangle, then all edges
of each four-cycle in $G$ belong to one triangle block.
\end{enumerate}
\end{proposition}

By Proposition~\ref{prop:blocks}(i), each copy of $K_4$ is a triangle
block, so distinct copies are edge-disjoint.
By Corollary~1.4 of \cite{CFSZ}, there is a set of
$o(n^{3/2})$ edges meeting every triangle.\footnote{The $o(n^{3/2})$
scale is essential here. The classical triangle-removal lemma
\cite{RS} would only give an $o(n^2)$ bound on the number of deleted
edges.} This set contains at least
two edges of each $K_4$, since deleting one edge leaves two triangles.
There are consequently $o(n^{3/2})$ copies of $K_4$. Deleting their
edges loses four triangles per copy, hence $o(n^{3/2})$ triangles.
Every surviving edge still lies in a triangle, and the remaining
blocks are books.

Fix a base edge in each of the books above, choosing any edge for
a one-triangle book. A \emph{retained book family} is obtained by
retaining an arbitrary subset of the triangles in each book and
discarding books with no remaining triangles. Each surviving book
keeps its original base, even when it has only one triangle.

Deleting a page from one book means removing its two spokes in
that book, together with the base edge if no triangle remains.
The page vertex remains, and none of its edges in other books
are removed. The union of a retained family is a subgraph of $G$,
and its triangles are exactly the retained book triangles.

Let $T$ be the number of retained triangles, and let the
\emph{triangle degree} $t_v$ be the number containing $v$. The
\emph{spoke graph} $H$ has the same vertex set as $G$ and consists
of the spokes of the retained books. Each spoke lies in exactly one
retained triangle, and each triangle contributes two spokes, so
$e(H)=2T$. The graph $H$ has no triangles, since every triangle in
the book union uses its base. It also has no five-cycles, since
it is a subgraph of $G$.

\subsection{The main lemmas}\label{subsec:main-lemmas}
The first lemma gives the upper bound when the spokes already form
a bipartite graph. By $e(H)=2T$, its leading constant for the
triangle count is $1/(3\sqrt3)$.

\begin{lemma}\label{lem:bipartite}
If the spoke graph of a retained book family on $n$ vertices is
bipartite, then
\[
 e(H)\le\left(\frac{2}{3\sqrt3}+o(1)\right)n^{3/2}.
\]
The error is uniform over all such families.
\end{lemma}

To find a partition with few spokes inside its parts, we first group
the vertices into classes and count the spokes between them. The
counts also record which endpoint of a spoke is a page. Orient every
spoke from its base endpoint to its page, and write $u\to v$ for an
edge with this orientation, which we call an \emph{arc}.

For a partition $V_1,\ldots,V_q$ of $V(H)$ into nonempty classes, set
$w_i=|V_i|/n$. Define the scaled adjacency and arc matrices by
\[
\begin{aligned}
 A_{ij}&=\frac{\sqrt n}{|V_i||V_j|}
   \bigl|\{(u,v)\in V_i\times V_j:uv\in E(H)\}\bigr|,\\
 X_{ij}&=\frac{\sqrt n}{|V_i||V_j|}
   \bigl|\{(u,v)\in V_i\times V_j:u\to v\}\bigr|.
\end{aligned}
\]

These counts use ordered pairs. In particular, a spoke inside $V_i$
is counted twice in the definition of $A_{ii}$ and once in that of
$X_{ii}$.

For a symmetric matrix $A$ on classes with weights $w_i$, define
\[
 W_5(A)=\sum_{i_1,\ldots,i_5}w_{i_1}\cdots w_{i_5}
 A_{i_1i_2}A_{i_2i_3}A_{i_3i_4}A_{i_4i_5}A_{i_5i_1}.
\]
When $A$ is nonnegative, this counts closed five-step walks with
weights, allowing repeated indices. Although the spoke graph has no
such walks, the weighted model may have them.

The next lemma gives a partition for which the total weight of the
closed five-step walks is small. In each row it also bounds the
weighted sum of $A_{ij}^2-X_{ji}^2/2$, with an error that is small
on average over the classes.
Write $\Delta(H)$ for the maximum degree of $H$.

\begin{lemma}\label{lem:model}
For every $C,\tau>0$, there are positive integers $K,n_0$,
depending only on $C,\tau$, such that every retained book family
on $n\ge n_0$ vertices with $\Delta(H)\le C\sqrt n$ has a
partition into at most $K$ nonempty classes whose matrices satisfy
the following.
\begin{enumerate}
\item[(i)] We have
\begin{equation}\label{eq:model-walk}
 W_5(A)<\tau.
\end{equation}
\item[(ii)] There are $\zeta_i\ge0$ such that, for every $i$,
\begin{equation}\label{eq:model-row}
 \sum_jw_j\bigl(A_{ij}^2-\tfrac12X_{ji}^2\bigr)\le1+\zeta_i,
 \quad \text{and}\quad \sum_iw_i\zeta_i<\tau.
\end{equation}
\end{enumerate}
\end{lemma}

For fixed $b>0$, if every vertex lies in at least $b\sqrt n$ triangles, then
$\Delta(H)=O_b(\sqrt n)$ by Lemma~\ref{lem:counts}(iv). We can then
use the partition from Lemma~\ref{lem:model} to prove the following
statement about the original vertices.

\begin{lemma}\label{lem:near}
Fix $c>1/(3\sqrt3)$ and $b>1/(2\sqrt3)$. For every $\eta>0$,
there is $n_0=n_0(b,c,\eta)$ such that every retained book family
on $n\ge n_0$ vertices satisfying
\[
 T\ge cn^{3/2},\qquad t_v\ge b\sqrt n\quad\text{for every }v
\]
admits a partition of its vertices into two parts such that fewer
than $\eta n^{3/2}$ spokes have both endpoints in the same part.
\end{lemma}

\subsection{The upper bound on the number of triangles}\label{sec:proof}
The construction of Bollob\'as and Gy\H{o}ri \cite{BG} gives the lower
bound. For the upper bound, we first consider a retained book family
on $n$ vertices with $T$ triangles.

\begin{proof}
Fix $c>1/(3\sqrt3)$. We shall choose a constant $K$ depending only
on $c$ and prove $T\le cn^{3/2}+K$ for every $n$. For the induction
step, we first show that, for all sufficiently large $n$,
\[
 T\ge cn^{3/2}\quad\Longrightarrow\quad
 \min_v t_v\le c\bigl(n^{3/2}-(n-1)^{3/2}\bigr).
\]
The right-hand side is the decrease in the proposed upper bound
when the number of vertices decreases from $n$ to $n-1$.

Suppose $T\ge cn^{3/2}$. Choose $b$ and $\eta$ in terms of $c$ so that
$1/(2\sqrt3)<b<3c/2$ and $0<\eta<c-1/(3\sqrt3)$. If $t_v\ge b\sqrt n$ for every
vertex $v$, Lemma~\ref{lem:near} gives a partition with fewer than
$\eta n^{3/2}$ spokes inside its parts. Discard every triangle
containing such a spoke. Fewer than $\eta n^{3/2}$ triangles are
discarded, since each spoke belongs to exactly one triangle.
All remaining spokes cross the partition, so their graph $H'$ is
bipartite on the same $n$ vertices. Each remaining triangle contributes
two spokes, so exactly $e(H')/2$ triangles remain. Hence
Lemma~\ref{lem:bipartite} gives
\[
 T\le \eta n^{3/2}+\tfrac12 e(H')
 \le \left(\eta+\frac1{3\sqrt3}+o(1)\right)n^{3/2}<cn^{3/2}
\]
for all sufficiently large $n$, a contradiction. Thus some vertex
satisfies $t_v<b\sqrt n$. Since
$c\bigl(n^{3/2}-(n-1)^{3/2}\bigr)/\sqrt n\to 3c/2>b$,
this vertex satisfies the claimed bound.

Fix $n_0\ge3$ depending only on $c$ so that the implication holds
for every $n\ge n_0$, and set $K=\binom{n_0}{3}$. We now prove
$T\le cn^{3/2}+K$ by induction on $n$. For $n<n_0$ this follows
from $T\le\binom n3\le K$. For
$n\ge n_0$ we may assume $T\ge cn^{3/2}$. Delete a vertex $v$ of
minimum triangle degree. This removes all triangles of each book in
which $v$ is a base endpoint and just one triangle from each book in
which $v$ is a page. The remaining triangles therefore form a retained
book family on $n-1$ vertices with exactly $T-t_v$ triangles.
By induction and the bound on $t_v$,
\[
 T\le c(n-1)^{3/2}+K+t_v\le cn^{3/2}+K.
\]
This completes the induction and gives $T\le cn^{3/2}+O_c(1)$.

The reduction in Section~\ref{sec:books} discards only $o(n^{3/2})$
triangles from the original $C_5$-free graph. Consequently
$\limsup_{n\to\infty}\ex(n,K_3,C_5)/n^{3/2}\le c$.
Letting $c\downarrow 1/(3\sqrt3)$ proves the upper bound.
\end{proof}

\section{Proof of Lemma~\ref{lem:bipartite}}\label{sec:bipartite}
We first count two-edge walks in an arbitrary retained book family,
keeping track of those with the same endpoints. These counts give
the neighborhood estimates used here and in the weighted approximation.
We then specialize to a bipartite spoke graph and prove the edge bound
in Lemma~\ref{lem:bipartite}.

\subsection{Neighborhood counts}\label{subsec:counts}
We work with the retained books of Section~\ref{sec:books}. The
following consequence of Proposition~\ref{prop:blocks} explains
how two of the original blocks can meet.

\begin{proposition}\label{prop:block-intersections}
Let $G$ be a $C_5$-free graph in which every edge lies in a triangle.
Then any two distinct triangle blocks of $G$ have at most one
common vertex.
\end{proposition}

\begin{proof}
Suppose that two distinct blocks share vertices $u,v$.
By Proposition~\ref{prop:blocks}(i), each block
contains a $u$--$v$ path of length two. The internal vertices of
these paths are distinct because the blocks are edge-disjoint.
Their union is therefore a four-cycle with edges in two distinct
blocks, contrary to Proposition~\ref{prop:blocks}(ii).
\end{proof}

For each book $F$ in the current family, write $\{s_F,t_F\}$ for
its base, $P_F$ for its set of remaining pages and $k_F=|P_F|$.

A four-cycle in the union of the retained books lies in one original
block by Proposition~\ref{prop:blocks}(ii). All four of its edges
therefore lie in the same retained book.

In the oriented spoke graph $H$, write $d(v),d^+(v),d^-(v)$ for the
degree, outdegree and indegree of $v$, and $N(v),N^-(v)$ for its
neighborhood and inneighborhood. We have $T=\sum_F k_F$ and
\begin{equation}\label{eq:roles}
 e(H)=2T,\qquad d(v)=d^+(v)+d^-(v),\qquad
 t_v=d^+(v)+\tfrac12d^-(v),\qquad \sum_vd(v)=4T.
\end{equation}
Each triangle contributes one to $d^+(v)$ at each base endpoint and two
to $d^-(v)$ at its page.

The spokes of each book $F$ form a $K_{2,k_F}$ with parts
$\{s_F,t_F\}$ and $P_F$. Suppose that distinct vertices $v,w$ lie
in the same part,
and choose a common neighbor $x$ in the opposite part. Any other
common neighbor $z$ in $H$ gives the four-cycle $v,x,w,z,v$,
which lies in a single retained book. This cycle contains the edge $vx$ of
$F$. Since that edge belongs
to only one retained book, the book containing the cycle is $F$.
Thus all common neighbors of $v,w$ lie in the opposite part of
this $K_{2,k_F}$. In particular, the base endpoints have exactly
$k_F$ common neighbors, even when $k_F=1$, and two pages have
exactly the two base endpoints as common neighbors.

Conversely, two distinct vertices with at least two common
neighbors lie on a four-cycle. Since this cycle lies in one book, the pair must
consist of its base endpoints or two of its pages. Every other
pair therefore has at most one common neighbor.

For the common-neighbor bound below, define a symmetric nonnegative
error matrix $E=(E(v,w))$. Its entries at coincident vertices and
at base pairs are
\[
 E(v,v)=d(v),\qquad
 E(s_F,t_F)=E(t_F,s_F)=k_F-1
\]
for every vertex $v$ and book $F$. Set all other entries to zero.

We next estimate the size of the second neighborhood. Let $N_2(u)$
be the set of vertices at distance exactly two from $u$ in $H$.
There are $\sum_{v\in N(u)}d(v)$ two-edge walks starting at $u$,
of which exactly $d(u)$ return to $u$. Since $H$ is triangle-free,
every other walk ends in $N_2(u)$, but several walks may have the
same endpoint.

If $u$ is a base endpoint of a book $F$, the $k_F$ walks through
its pages all reach the other base endpoint, giving $k_F-1$
repetitions. If $u$ is a page of $F$, each of the other $k_F-1$
pages is reached through both base endpoints, again giving
$k_F-1$ repetitions. By the common-neighbor classification above, every repeated
endpoint arises in one of these two ways. Therefore
\begin{equation}\label{eq:qexact}
 |N_2(u)|=\sum_{v\in N(u)}d(v)-d(u)
       -\sum_{F:\,u\in\{s_F,t_F\}\cup P_F}(k_F-1).
\end{equation}

The repetitions at base endpoints satisfy
\[
 \sum_{F:\,u\in\{s_F,t_F\}}(k_F-1)
 \le\sum_{F:\,u\in\{s_F,t_F\}}k_F=d^+(u).
\]
When $u$ is a page of $F$, both base endpoints have outdegree at
least $k_F$. Summing over these books gives
\[
 \sum_{F:\,u\in P_F}(k_F-1)
 \le\tfrac12\sum_{F:\,u\in P_F}
             \bigl(d^+(s_F)+d^+(t_F)\bigr)
 =\tfrac12\sum_{v\to u}d^+(v).
\]
The last equality holds because each arc into $u$ joins a base
endpoint to the page $u$ in a unique book.
Substituting these bounds into \eqref{eq:qexact}, we obtain
\[
 |N_2(u)|\ge
 \left(\sum_{v\in N(u)}d(v)-\tfrac12\sum_{v\to u}d^+(v)\right)
 -d(u)-d^+(u).
\]
The return walks and the repetitions at base endpoints together
contribute at most $d(u)+d^+(u)\le2d(u)$.
For the expression in parentheses, write
\begin{equation}\label{eq:h}
 h(u)=\sum_{v\in N(u)}d(v)-\tfrac12\sum_{v\to u}d^+(v).
\end{equation}
Thus $|N_2(u)|\ge h(u)-2d(u)$.

The following lemma collects this second-neighborhood estimate together with
the degree and common-neighbor bounds.

\begin{lemma}\label{lem:counts}
For every retained book family on $n$ vertices, the following hold.
\begin{enumerate}
\item[(i)] We have
\begin{equation}\label{eq:moment}
 \sum_vd(v)^2=O(n^2),\qquad T=O(n^{3/2}).
\end{equation}
\item[(ii)] For all vertices $v,w$ (including $v=w$),
\begin{equation}\label{eq:cap}
 |N(v)\cap N(w)|
 -\tfrac12|N^-(v)\cap N^-(w)|\le1+E(v,w),
\end{equation}
and these errors satisfy
\begin{equation}\label{eq:error-sums}
 \sum_{v,w}E(v,w)\le6T,\qquad
 \sum_wE(v,w)\le2d(v)\quad\text{for every }v.
\end{equation}
\item[(iii)] For every vertex $v$,
\begin{equation}\label{eq:q}
 |N_2(v)|\ge h(v)-2d(v),\qquad
 h(v)\ge\tfrac12\sum_{u\in N(v)}d(u).
\end{equation}
\item[(iv)] For each fixed $b>0$, if $t_v\ge b\sqrt n$ at every
vertex, then $\Delta(H)=O_b(\sqrt n)$.
\end{enumerate}
\end{lemma}
\begin{proof}
The sum $\sum_vd(v)^2$ counts ordered two-edge walks, including those
whose endpoints coincide. For distinct endpoints, the common-neighbor
classification gives an upper bound of two, except at base pairs.
For the ordered base pairs, adding $2\sum_Fk_F=2T$ to this bound
is enough. Walks with coincident endpoints contribute
$\sum_vd(v)=4T$. Hence
\[
 \sum_vd(v)^2\le2n(n-1)+6T.
\]
Since $2T=e(H)\le n^2/2$, the first assertion of \eqref{eq:moment}
follows. Cauchy--Schwarz and \eqref{eq:roles} give the bound on $T$.

We next prove the common-neighbor estimate. For two distinct pages
of one book, both common neighbors point to both pages,
so the left side of \eqref{eq:cap} is $2-2/2=1$. All other distinct
nonbase pairs have at most one common neighbor. For a base pair,
the left side of \eqref{eq:cap} is at most $k_F=1+E(s_F,t_F)$,
and for $v=w$ it is at most $d(v)=E(v,v)$.
The base pairs are distinct, so their total contribution is at most
$2T$, while the diagonal contributes $4T$. Moreover,
\[
 \sum_wE(v,w)
 =d(v)+\sum_{F:\,v\in\{s_F,t_F\}}(k_F-1)
 \le d(v)+d^+(v)\le2d(v).
\]
This proves \eqref{eq:cap}--\eqref{eq:error-sums}.

The first inequality in \eqref{eq:q} was established above.
The second follows from \eqref{eq:h} and
$\sum_{u\to v}d^+(u)\le\sum_{u\in N(v)}d(u)$.

Finally, if $t_u\ge b\sqrt n$ for every $u$, then $d(u)\ge b\sqrt n$.
Consequently \eqref{eq:q} gives
\[
 n\ge |N_2(v)|\ge(b\sqrt n/2-2)d(v).
\]
For sufficiently large $n$, this implies $d(v)=O_b(\sqrt n)$.
\end{proof}

\subsection{Counting spokes in the bipartite case}
We use H\"older's inequality to bound a sum of degree products from
below, and the second-neighborhood estimates to bound the same sum
from above. Comparing the bounds will give Lemma~\ref{lem:bipartite}.

\begin{proof}[Proof of Lemma~\ref{lem:bipartite}]
Fix a bipartition $V(H)=A\sqcup B$, and put
$a=|A|$, $b=|B|$ and $m=e(H)$. We may assume that $m>0$, so
$a,b>0$. We bound $m$ by estimating the same quantity from below
and above. Set
\[
 f(v)=d^+(v)+\tfrac23d^-(v),\qquad
 Q=\tfrac32\sum_{uv\in E(H)}f(u)f(v),
\]
where each edge is counted once.

For the lower bound, separate arcs starting in $A$ from arcs
starting in $B$. Each summand below equals one, and all denominators are
positive, so
\[
\begin{aligned}
 m={}\sum_{\substack{s\to t\\s\in A}}
   \left(\tfrac32f(s)f(t)\right)^{1/3}
   \left(\frac1{f(s)}\right)^{1/3}
   \left(\frac2{3f(t)}\right)^{1/3}
 +\sum_{\substack{s\to t\\s\in B}}
   \left(\tfrac32f(s)f(t)\right)^{1/3}
   \left(\frac2{3f(t)}\right)^{1/3}
   \left(\frac1{f(s)}\right)^{1/3}.
\end{aligned}
\]
In each summand, the second and third factors correspond to the
endpoints in $A$ and $B$, respectively. Applying H\"older's inequality
once to the full set of arcs gives
\[
\begin{aligned}
 m^3\le Q
 \left(
   \sum_{\substack{s\to t\\s\in A}}\frac1{f(s)}
   +\sum_{\substack{s\to t\\s\in B}}\frac2{3f(t)}
 \right)
 \cdot
 \left(
   \sum_{\substack{s\to t\\s\in A}}\frac2{3f(t)}
   +\sum_{\substack{s\to t\\s\in B}}\frac1{f(s)}
 \right).
\end{aligned}
\]
Grouping the second factor by its endpoint in $A$, we obtain
\[
\begin{aligned}
 \sum_{\substack{s\to t\\s\in A}}\frac1{f(s)}
  +\sum_{\substack{s\to t\\s\in B}}\frac2{3f(t)}
 =
 \sum_{\substack{u\in A\\f(u)>0}}
   \frac{d^+(u)+\tfrac23d^-(u)}{f(u)}
 =\sum_{\substack{u\in A\\f(u)>0}}1\le a.
\end{aligned}
\]
Similarly, grouping the third factor by its endpoint in $B$ gives
\[
\begin{aligned}
 \sum_{\substack{s\to t\\s\in A}}\frac2{3f(t)}
  +\sum_{\substack{s\to t\\s\in B}}\frac1{f(s)}
=
 \sum_{\substack{v\in B\\f(v)>0}}
   \frac{\tfrac23d^-(v)+d^+(v)}{f(v)}
 =\sum_{\substack{v\in B\\f(v)>0}}1\le b.
\end{aligned}
\]
Consequently,
\begin{equation}\label{eq:product}
 m^3
 \le abQ.
\end{equation}

For the upper bound, let $x$ be the number of arcs from $A$ to $B$.
We will show that $Q\le2nx/3+b(m-x)+O(n^2)$ using second-neighborhood
estimates. The base endpoints of each book lie in the same part,
since they share a page neighbor, and its pages lie in the other part.
Consider a book $F$ with base endpoints $s,t\in A$. Both $N_2(s)$
and $N_2(t)$ lie in $A$, and we claim that
\begin{equation}\label{eq:overlap}
 |N_2(s)|+|N_2(t)|\le a+\sum_{p\in P_F}d(p).
\end{equation}
To bound their intersection, take $z\in N_2(s)\cap N_2(t)$ and
choose paths $spz$ and $tp'z$. If $p\ne p'$, these paths together
with the base edge $ts$ form the five-cycle $s,p,z,p',t,s$ in the
book union. Its vertices are distinct because
$z\in A\setminus\{s,t\}$ and $p,p'\in B$.
Thus $p=p'$, and the common-neighbor classification shows that
$p\in P_F$. Every vertex of the intersection is therefore adjacent
to a page of $F$, so the intersection has size at most
$\sum_{p\in P_F}d(p)$. Since the union is contained in $A$,
\eqref{eq:overlap} follows.

Multiply \eqref{eq:overlap} by $2k_F$ and sum over books based in
$A$. The coefficient of $|N_2(u)|$ is $2d^+(u)$, while
$2\sum_{F:\{s_F,t_F\}\subseteq A}k_F=x$. Since $d^+(s_F)+d^+(t_F)\ge2k_F$, the term $2k_Fd(p)$ for each
page $p$ is at most $d^+(s_F)d(p)+d^+(t_F)d(p)$.
Summing over pages and books gives
\[
 2\sum_{u\in A}d^+(u)|N_2(u)|
 -\sum_{\substack{s\to t\\s\in A}}d^+(s)d(t)\le ax.
\]
By \eqref{eq:q}, we may replace $|N_2(u)|$ by $h(u)$ if we add
at most $4\sum_ud^+(u)d(u)=O(n^2)$ to the right-hand side, using
\eqref{eq:moment}.
Exchanging the parts gives a second estimate. We also use
$|N_2(v)|\le b$ for $v\in B$, together with
$\sum_{v\in B}d^-(v)=x$. Thus
\begin{equation}\label{eq:budgets}
\begin{aligned}
 2\sum_{u\in A}d^+(u)h(u)
 -\sum_{\substack{s\to t\\s\in A}}d^+(s)d(t)
 &\le ax+O(n^2),\\
 2\sum_{v\in B}d^+(v)h(v)
 -\sum_{\substack{s\to t\\s\in B}}d^+(s)d(t)
 &\le b(m-x)+O(n^2),\\
 \sum_{v\in B}d^-(v)h(v)&\le bx+O(n^2).
\end{aligned}
\end{equation}

Let $L$ be the sum of the three left-hand sides in
\eqref{eq:budgets}, with coefficients $2/3,1,2/3$, respectively.
The first and third right-hand sides combine into $2nx/3$, so
\[
 L\le\tfrac23nx+b(m-x)+O(n^2).
\]
It remains to check that $Q\le L$. We do this one arc at a time.
By \eqref{eq:h}, an arc $s\to t$ contributes $d(t)$ to $h(s)$
and $d(s)-d^+(s)/2$ to $h(t)$. Its contribution to $L$, including
the subtracted arc sums, is therefore
\[
\begin{cases}
 \tfrac23d^+(s)d(t)
 +(f(t)+d^+(t))(d(s)-\tfrac12d^+(s)),&s\in A,\\[2pt]
 f(s)d(t)+\tfrac43d^+(t)(d(s)-\tfrac12d^+(s)),&s\in B.
\end{cases}
\]
Substituting $d(v)=3f(v)/2-d^+(v)/2$ into either expression gives
the right-hand side below. Its last term is nonnegative, so
\[
 \tfrac32f(s)f(t)\le\tfrac32f(s)f(t)
 +d^+(t)\bigl(d^-(s)+\tfrac16d^+(s)\bigr).
\]
Summing over arcs gives $Q\le L$, as required. Since $Q$ is
unchanged when the parts are exchanged, the same argument gives
$Q\le ax+2n(m-x)/3+O(n^2)$. Together with \eqref{eq:product},
we obtain
\begin{equation}\label{eq:combined-bound}
 \frac{m^3}{ab}\le Q\le
 \min\{\tfrac23nx+b(m-x),\ ax+\tfrac23n(m-x)\}+O(n^2).
\end{equation}

We next bound the right-hand side using $a+b=n$, so that the
resulting estimate depends only on $n$ and $m$.
If $a\le n/3$, then $b\ge2n/3$, so the first expression in the
minimum is at most $bm$. After multiplication by $ab$, it is at
most $ab^2m\le4n^3m/27$, by the arithmetic--geometric mean
inequality applied to $a,b/2,b/2$. The case $b\le n/3$ is symmetric.
Otherwise $a,b\ge n/3$. Give the first expression weight
$(2b-a)/n$ and the second weight $(2a-b)/n$. These nonnegative
weights sum to one and cancel $x$, giving the weighted average
$m(4n/3-3ab/n)$.
Thus in this case the minimum, multiplied by $ab$, is at most
\[
 m\,ab\left(\frac{4n}3-\frac{3ab}n\right)
 =m\left(\frac4{27}n^3-\frac3n
                  \left(ab-\frac29n^2\right)^2\right)
 \le\frac4{27}n^3m.
\]
In all cases, \eqref{eq:combined-bound} gives
$m^3\le4n^3m/27+O(n^4)$.
If $m<2n^{3/2}/(3\sqrt3)$, the desired bound already holds.
Otherwise $m=\Omega(n^{3/2})$, and division by $m$ gives
$m^2\le4n^3/27+O(n^{5/2})$, as required.
\end{proof}

\section{Proof of Lemma~\ref{lem:model}}\label{sec:matrix}
The matrices $A,X$ from Section~\ref{subsec:main-lemmas} record the
total number of triangles and the average triangle degree in each
class. We first verify these identities for an arbitrary partition.
We then prove Lemma~\ref{lem:model} by choosing a partition that also
satisfies the bounds on closed five-step walks and on the rows.

For $i\ne j$, every edge between $V_i$ and $V_j$ has exactly
one of the two orientations, so $A_{ij}=X_{ij}+X_{ji}$.
The diagonal entries record edges with both endpoints in the
same class. Each such edge gives two ordered pairs in the count
for $A_{ii}$, but only one arc in the count for $X_{ii}$. Hence
\[
 A_{ii}=\frac{2\sqrt n\,e(H[V_i])}{|V_i|^2},\qquad
 X_{ii}=\frac{\sqrt n\,e(H[V_i])}{|V_i|^2},
\]
where $e(H[V_i])$ is the number of edges with both endpoints in
$V_i$. Thus $A_{ii}=2X_{ii}$, and $A=X+X^{\mathsf T}$ holds on
the diagonal as well.

Since $e(H)=2T$ and $t_v=d^+(v)+d^-(v)/2$, every partition satisfies
\begin{equation}\label{eq:model-identities}
 \sum_{i,j}w_iw_jA_{ij}=\frac{4T}{n^{3/2}},\qquad
 \sum_jw_j(X_{ij}+X_{ji}/2)
 =\frac1{|V_i|\sqrt n}\sum_{v\in V_i}t_v.
\end{equation}
Thus, if $t_v\ge b\sqrt n$ at every vertex, then each row sum
in \eqref{eq:model-identities} is at least $b$.

We use weighted matrix operations on classes of unequal sizes.
Here $f,g\colon\{1,\ldots,q\}\to\mathbb R$ are functions on the
class indices, and $Y,Z\in\mathbb R^{q\times q}$ are matrices
indexed by the same classes.
For positive weights $w_i$ summing to one, set
$\langle f,g\rangle=\sum_iw_if(i)g(i)$, with
$\|f\|_2^2=\langle f,f\rangle$, and let
$(Yf)(i)=\sum_jw_jY_{ij}f(j)$. Thus
\[
 (YZ)_{ij}=\sum_kw_kY_{ik}Z_{kj},\qquad
 \|Y\|_{\Frob}^2=\sum_{i,j}w_iw_jY_{ij}^2,
 \qquad \tr Y=\sum_iw_iY_{ii}.
\]
Here $\|Y\|_{\Frob}$ is the weighted Frobenius norm.
The operator norm is $\|Y\|_{\op}=\sup_{\|f\|_2=1}\|Yf\|_2$.
These are the usual operations and norms of the representing
matrix $(\sqrt{w_iw_j}Y_{ij})$. For symmetric $Y$, the operator
norm is the largest absolute eigenvalue and $\tr(Y^4)$ is the
sum of the fourth powers of its eigenvalues.
In particular, $W_5(A)=\tr(A^5)$.

\begin{proof}[Proof of Lemma~\ref{lem:model}]
We may assume that $0<\tau<1$. A partition $\mathcal Q$ is a
\emph{refinement} of $\mathcal P$ if every class of $\mathcal Q$
is contained in a class of $\mathcal P$.
We first find a partition such that every refinement still has
a small five-step walk count.
We then choose a refinement satisfying \eqref{eq:model-row}.

All norms on the original vertex set use weights $1/n$.
For a partition $\mathcal P=\{V_1,\ldots,V_q\}$ of $V(H)$,
define $P$ on functions $f\colon V(H)\to\mathbb R$ by
\[
 (Pf)(v)=\frac1{|V_i|}\sum_{u\in V_i}f(u)
 \qquad\text{when }v\in V_i.
\]
The function $Pf$ is constant on each class, while $f-Pf$ has
mean zero on each class. Thus $P$ is the orthogonal projection
onto the space of functions constant on each class.

\medskip\noindent
\emph{Claim 1.} There is a partition $\mathcal P_0$ with a number of
classes bounded in terms of $C,\tau$ such that every refinement
of $\mathcal P_0$ satisfies $W_5(A)<\tau/2$.

\begin{proof}[Proof of Claim~1]
Put $A_n(u,v)=\sqrt n\,\mathbf1_{uv\in E(H)}$.
Matrix multiplication uses vertex weights $1/n$, so
\[
 (A_n^2)_{uv}=\frac1n\sum_xA_n(u,x)A_n(x,v)
 =|N(u)\cap N(v)|.
\]
Fix $u$. If $v\ne u$ and $\{u,v\}$ is not a book base, then
$u,v$ have at most two common neighbors. Hence each squared entry
$(A_n^2)_{uv}^2$ is at most $4$, and all such vertices together
contribute at most $4n$.
If $\{u,v\}$ is the base of a book $F$, its contribution is $k_F^2$.
These contributions sum to at most
$\sum_{F:\,u\in\{s_F,t_F\}}k_F^2\le(d^+(u))^2$,
since $d^+(u)$ is the sum of the page counts of all books with
$u$ as a base endpoint. The diagonal contributes $d(u)^2$.
Using $\Delta(H)\le C\sqrt n$ and the symmetry of $A_n$, we obtain
\begin{equation}\label{eq:row-four}
\begin{aligned}
 \frac1n\sum_v(A_n^2)_{uv}^2
 &\le4+\frac{(d^+(u))^2+d(u)^2}{n}\le M:=4+2C^2,\\
 \tr(A_n^4)&=\frac1n\sum_u(A_n^4)_{uu}
 =\frac1{n^2}\sum_{u,v}(A_n^2)_{uv}^2\le M.
\end{aligned}
\end{equation}
Here $(A_n^4)_{uu}=\frac1n\sum_v(A_n^2)_{uv}(A_n^2)_{vu}$,
and $A_n^2$ is symmetric. The trace bound follows by averaging
the first inequality over $u$.

We construct classes on which the eigenvectors with large
absolute eigenvalues vary little. Averaging then changes the terms
belonging to these eigenvectors
only slightly. The sum of the remaining terms already has small
operator norm.
Choose an orthonormal eigenbasis $f_1,\ldots,f_n$ of $A_n$, with
eigenvalues $\lambda_1,\ldots,\lambda_n$, and retain the vectors
whose eigenvalues satisfy $|\lambda_r|>\gamma$, for a parameter
$\gamma>0$. By the trace bound in \eqref{eq:row-four}, there are
at most $M/\gamma^4$ retained vectors.
The row bound also bounds the absolute values of their entries.
For a retained unit eigenvector $f$ with eigenvalue $\lambda$,
we have $A_n^2f=\lambda^2f$ and $\frac1n\sum_v f(v)^2=1$.
Thus Cauchy--Schwarz gives
\[
\begin{aligned}
 \lambda^2|f(u)|
=\left|\frac1n\sum_v(A_n^2)_{uv}f(v)\right|
\le\left(\frac1n\sum_v(A_n^2)_{uv}^2\right)^{1/2}
      \left(\frac1n\sum_v f(v)^2\right)^{1/2}
 \le\sqrt M.
\end{aligned}
\]
Thus these values lie in
$[-\sqrt M/\gamma^2,\sqrt M/\gamma^2]$.
For $\ell>0$, divide this interval into disjoint intervals of
length at most $\ell$. Take the intervals to be left-closed and
right-open, except that the last one includes the right endpoint.
At most $\lceil2\sqrt M/(\gamma^2\ell)\rceil$ intervals are needed.
Put two vertices in the same class when their values lie in the
same interval for every retained vector. The resulting partition
$\mathcal P_0$ has a class count bounded in terms of
$C,\gamma,\ell$, since there are at most $M/\gamma^4$ retained
vectors. If none are retained, take $\mathcal P_0=\{V(H)\}$.

For the averaging projection $P$ of any refinement of
$\mathcal P_0$, each retained vector still satisfies
$\|f-Pf\|_2\le\ell$, since splitting a class cannot enlarge the
range of its values.
For a unit function $g=\sum_{r=1}^n a_rf_r$, we have
$\sum_r a_r^2=1$ and
$(A_n-PA_n)g=\sum_r\lambda_ra_r(f_r-Pf_r)$.
Since $P$ is an orthogonal projection,
$\|h-Ph\|_2\le\|h\|_2$ for every function $h$.
Splitting the eigenvector sum at $\gamma$ and using the triangle
inequality therefore gives
\[
\begin{aligned}
 \|(A_n-PA_n)g\|_2
 &\le\left\|\sum_{|\lambda_r|\le\gamma}\lambda_ra_rf_r\right\|_2
       +\ell\sum_{|\lambda_r|>\gamma}|\lambda_r||a_r|
 =\left(\sum_{|\lambda_r|\le\gamma}\lambda_r^2a_r^2\right)^{1/2}
       +\ell\sum_{|\lambda_r|>\gamma}|\lambda_r||a_r|\\
 &\le\gamma+\ell\sum_{|\lambda_r|>\gamma}|\lambda_r|.
\end{aligned}
\]
The equality uses orthonormality, and the last inequality uses
$\sum_r a_r^2=1$ and $|a_r|\le1$.
Taking the supremum over unit $g$ and using
$\sum_r|\lambda_r|^4=\tr(A_n^4)\le M$, we obtain
\begin{equation}\label{eq:approx}
 \|A_n-PA_n\|_{\op}
 \le\gamma+\ell\sum_{|\lambda_r|>\gamma}|\lambda_r|
 \le\gamma+\frac{M\ell}{\gamma^3}.
\end{equation}
The last inequality uses
$|\lambda_r|\le|\lambda_r|^4/\gamma^3$ for the retained eigenvalues.

We now bound $W_5(A)$ by $5M\|A_n-PA_n\|_{\op}$.
For $u\in V_i$ and $v\in V_j$, averaging over both classes gives
\[
 (PA_nP)_{uv}
 =\frac1{|V_i||V_j|}\sum_{x\in V_i,\,y\in V_j}A_n(x,y)
 =A_{ij}.
\]
Grouping the vertex sums in the trace by their classes gives
the first equality below. The second uses $P^2=P$ and cyclicity
of trace.
\[
 W_5(A)=\tr((PA_nP)^5)=\tr((PA_n)^5).
\]
The graph $H$ has neither triangles nor five-cycles, so it has no
closed five-step walks and $\tr(A_n^5)=0$.

We compare these traces using the cycle counting lemma,
Lemma~51 of \cite{Pascadi}.
For a real matrix $Y$, its Schatten $4$-norm is
$\|Y\|_{S_4}=(\tr((YY^{\mathsf T})^2))^{1/4}$,
with the weighted multiplication and trace defined above.
Averaging the first variable does not increase this norm by
Lemma~37(iii) of \cite{Pascadi}, so \eqref{eq:row-four} gives
\[
 \|PA_n\|_{S_4}\le\|A_n\|_{S_4}
 =\tr(A_n^4)^{1/4}\le M^{1/4}.
\]
Apply Lemma~51 on five copies of $V(H)$ with uniform vertex
weights $1/n$, comparing $A_n$ with $PA_n$ in each factor.
Take $k=5$, $q_i=4$ and $r_i=\infty$ for $1\le i\le5$, with
norm bound $M^{1/4}$. The Schatten $\infty$-norm is the operator
norm, and the required exponent condition is $0+4/4=1$.
Since $W_5(A)\ge0$ and $\tr(A_n^5)=0$, the lemma yields
\begin{equation}\label{eq:five}
 0\le W_5(A)
 =\bigl|\tr((PA_n)^5)-\tr(A_n^5)\bigr|
 \le5M\|A_n-PA_n\|_{\op}.
\end{equation}
Choose $0<\gamma<\tau/(20M)$ and
$0<\ell<\tau\gamma^3/(20M^2)$, both depending only on $C,\tau$.
Then $\gamma+M\ell/\gamma^3<\tau/(10M)$, so
\eqref{eq:approx} and \eqref{eq:five} give $W_5(A)<\tau/2$
for every refinement of $\mathcal P_0$. The number of classes in
$\mathcal P_0$ is therefore bounded in terms of $C,\tau$ alone.
\end{proof}

\medskip\noindent
\emph{Claim 2.} Every partition into $q_0$ classes has a refinement
satisfying \eqref{eq:model-row}, with a number of classes bounded
in terms of $C,\tau,q_0$, provided $n$ is sufficiently large in
terms of these parameters.

\begin{proof}[Proof of Claim~2]
Taking class averages decreases both squared terms in the row
inequality. The decrease in the term with a positive sign helps
the upper bound, but the decrease in the subtracted term must be
included in the error.
We therefore consider the outdegrees into each class, whose
squared averages form the subtracted term. For a current partition
$\mathcal P=\{V_1,\ldots,V_q\}$,
write $d(v,S)=|N(v)\cap S|$ and
$d^+(v,S)=|\{u\in S:v\to u\}|$, and put
\[
 g_i(v)=\frac{\sqrt n\,d^+(v,V_i)}{|V_i|}.
\]
Since $X_{ji}$ counts arcs from $V_j$ to $V_i$, we have
\[
 \frac1{|V_j|}\sum_{v\in V_j}g_i(v)
 =\frac{\sqrt n}{|V_i||V_j|}\sum_{v\in V_j}d^+(v,V_i)
 =X_{ji}.
\]
For fixed $v$, the square $d^+(v,V_i)^2$ counts ordered pairs
$(s,t)\in V_i\times V_i$ with $v\to s$ and $v\to t$, allowing
$s=t$. For fixed $s,t$, the possible vertices $v$ are their common
in-neighbors, forming $N^-(s)\cap N^-(t)$.
Counting these ordered triples in the two ways gives
\[
\begin{aligned}
 \|g_i\|_2^2
 =\frac1n\sum_v g_i(v)^2
 =\frac1{|V_i|^2}\sum_v d^+(v,V_i)^2
 &=\frac1{|V_i|^2}\sum_{s,t\in V_i}|N^-(s)\cap N^-(t)|.
\end{aligned}
\]
For the undirected count, $A_{ij}$ is the average of
$\sqrt n\,d(v,V_i)/|V_i|$ over $v\in V_j$.
The square of an average is at most the average of the squares, so
\[
\begin{aligned}
 \sum_jw_jA_{ij}^2
 &=\sum_j\frac{|V_j|}{n}
   \left(\frac1{|V_j|}\sum_{v\in V_j}
          \frac{\sqrt n\,d(v,V_i)}{|V_i|}\right)^2\\
 &\le\frac1n\sum_v
          \left(\frac{\sqrt n\,d(v,V_i)}{|V_i|}\right)^2
 =\frac1{|V_i|^2}\sum_{s,t\in V_i}|N(s)\cap N(t)|.
\end{aligned}
\]
The last equality counts ordered triples $(v,s,t)$ with
$s,t\in V_i$ and $vs,vt\in E(H)$, again allowing $s=t$.
Subtracting half the identity for $\|g_i\|_2^2$ and applying
\eqref{eq:cap} to each ordered pair $(s,t)$ gives
\[
\begin{aligned}
 \sum_jw_jA_{ij}^2-\tfrac12\|g_i\|_2^2
 \le\frac1{|V_i|^2}\sum_{s,t\in V_i}
       \bigl(|N(s)\cap N(t)|-\tfrac12|N^-(s)\cap N^-(t)|\bigr)
 \le1+\frac1{|V_i|^2}\sum_{s,t\in V_i}E(s,t).
\end{aligned}
\]
To replace $\|g_i\|_2^2$ by the squared class averages, note that
$Pg_i=X_{ji}$ on $V_j$ and
$\sum_{v\in V_j}(g_i(v)-X_{ji})=0$.
Expanding $g_i(v)^2=(X_{ji}+g_i(v)-X_{ji})^2$ on each class
therefore gives
\[
\begin{aligned}
 \|g_i\|_2^2
 =\sum_jw_jX_{ji}^2
   +\frac1n\sum_j\sum_{v\in V_j}(g_i(v)-X_{ji})^2
 =\sum_jw_jX_{ji}^2+\|g_i-Pg_i\|_2^2.
\end{aligned}
\]
Substituting this identity into the preceding inequality gives
the row inequality in \eqref{eq:model-row} with
\[
 \zeta_i=\frac1{|V_i|^2}\sum_{s,t\in V_i}E(s,t)
          +\tfrac12\|g_i-Pg_i\|_2^2.
\]
The first term averages the errors $E(s,t)$ in the common-neighbor
bound over $V_i\times V_i$. The second is half the mean squared
difference between $g_i$ and its class averages.

By \eqref{eq:error-sums} and $d(s)\le C\sqrt n$, the first terms satisfy
\[
\begin{aligned}
 \sum_i\frac{w_i}{|V_i|^2}\sum_{s,t\in V_i}E(s,t)
 \le\frac1n\sum_i\frac1{|V_i|}\sum_{s\in V_i}2d(s)
 \le\frac1n\sum_i2C\sqrt n=\frac{2Cq}{\sqrt n}.
\end{aligned}
\]
Summing the two terms in $\zeta_i$, we obtain
\begin{equation}\label{eq:averaged-error}
 \sum_iw_i\zeta_i\le\frac{2Cq}{\sqrt n}
       +\tfrac12\sum_iw_i\|g_i-Pg_i\|_2^2.
\end{equation}
It remains to find a refinement with $q$ bounded in terms of
$C,\tau,q_0$ and $\sum_iw_i\|g_i-Pg_i\|_2^2\le\tau$.
The bound on $q$ will then allow us to choose $n_0$ so that
$2Cq/\sqrt n<\tau/2$ for every $n\ge n_0$.

The functions $g_i$ count arcs into the current classes, so they
must be redefined after each refinement. Making the old functions
close to their new class averages therefore does not yet give the
required error bound for the new functions. To ensure that the
refinement process stops, we use the energy
\[
 \mathcal E(\mathcal P)=\sum_{i,j}w_iw_jX_{ij}^2
                      =\sum_iw_i\|Pg_i\|_2^2.
\]
Whenever the required error bound fails, we will construct a
refinement that increases this energy by more than $\tau/2$.

To bound this energy, recall that two distinct vertices have at
most two common predecessors.
Indeed, a base pair has none, two pages of a book have exactly
two, and every other pair has at most one common neighbor.
Together with $d^-(s)\le C\sqrt n$, this gives
\[
\begin{aligned}
 w_i\|g_i\|_2^2
 =\frac1{n|V_i|}\sum_{s,t\in V_i}|N^-(s)\cap N^-(t)|
 \le\frac{2|V_i|(|V_i|-1)+\sum_{s\in V_i}d^-(s)}{n|V_i|}
 \le2w_i+\frac C{\sqrt n}.
\end{aligned}
\]
Here the diagonal terms $s=t$ contribute $d^-(s)$, and each of
the $|V_i|(|V_i|-1)$ other ordered pairs contributes at most two.
Consequently,
\begin{equation}\label{eq:finite-energy}
 \mathcal E(\mathcal P)\le\sum_iw_i\|g_i\|_2^2
                       \le2+\frac{Cq}{\sqrt n}.
\end{equation}

We now construct one refinement, keeping the current functions
$g_i$ fixed. Group vertices with similar values of each $g_i$,
and set $\varepsilon=\tau/(2Cq)$.
Since $d^+(v)\le C\sqrt n$, we have
$0\le g_i(v)\le C/w_i$.
Divide $[0,C/w_i]$ into intervals of length $\varepsilon/w_i$,
shortening the last one if necessary. Use left-closed, right-open
intervals, with the final right endpoint included. Split each old class
according to the intervals containing $g_1(v),\ldots,g_q(v)$.
This choice of interval length gives at most
$1+\lceil C/\varepsilon\rceil$ intervals for each function,
independently of the class size. The resulting partition
$\mathcal Q$ has at most $q(1+\lceil C/\varepsilon\rceil)^q$ classes.

Let $Q$ be the averaging projection of $\mathcal Q$, and let
$s_i(v)$ be the lower endpoint of the interval containing
$g_i(v)$. It is constant on each new class, and
$0\le g_i-s_i\le\varepsilon/w_i$. Since $s_i\ge0$,
$(g_i-s_i)^2\le(\varepsilon/w_i)g_i$.
Among functions that are constant on each new class, $Qg_i$
minimizes the mean squared difference from $g_i$. Indeed,
$g_i-Qg_i$ has mean zero on each new class, so
\[
 \|g_i-s_i\|_2^2
 =\|g_i-Qg_i\|_2^2+\|Qg_i-s_i\|_2^2
 \ge\|g_i-Qg_i\|_2^2.
\]
It follows that
\[
\begin{aligned}
 w_i\|g_i-Qg_i\|_2^2
 &\le w_i\|g_i-s_i\|_2^2
 \le\frac{\varepsilon}{n}\sum_vg_i(v)
 =\frac{\varepsilon}{|V_i|\sqrt n}
      \sum_{u\in V_i}d^-(u)
 \le C\varepsilon.
\end{aligned}
\]
Summing over $i$ gives
$\sum_iw_i\|g_i-Qg_i\|_2^2\le Cq\varepsilon=\tau/2$.

We next show that this refinement gives the claimed increase in
energy. Each old class $V_i$ is a disjoint union of new classes,
so
\[
 g_i(v)=\sum_{\substack{W\in\mathcal Q\\W\subseteq V_i}}
       \frac{|W|}{|V_i|}\,\frac{\sqrt n\,d^+(v,W)}{|W|}.
\]
Averaging over the new classes preserves this convex combination.
Applying Jensen's inequality and using
$w_i|W|/|V_i|=|W|/n$, we obtain
\[
 \sum_iw_i\|Qg_i\|_2^2
 \le\sum_{W\in\mathcal Q}\frac{|W|}{n}
     \left\|Q\left(\frac{\sqrt n\,d^+(\cdot,W)}{|W|}\right)\right\|_2^2
 =\mathcal E(\mathcal Q).
\]
The last equality is the energy identity for the new target
classes $W$. Thus splitting the target classes can only increase
the sum of squared averages.
Using the squared-norm identity above for both $P$ and $Q$ now gives
\[
\begin{aligned}
 \mathcal E(\mathcal Q)-\mathcal E(\mathcal P)
 \ge\sum_iw_i\bigl(\|Qg_i\|_2^2-\|Pg_i\|_2^2\bigr)
 =\sum_iw_i\|g_i-Pg_i\|_2^2
       -\sum_iw_i\|g_i-Qg_i\|_2^2
 >\tau/2
\end{aligned}
\]
whenever $\sum_iw_i\|g_i-Pg_i\|_2^2>\tau$.

Starting from the given partition, stop when
$\sum_iw_i\|g_i-Pg_i\|_2^2\le\tau$; otherwise apply the refinement
above and define the functions $g_i$ from the new classes.
To prove that this process stops, consider its first
$\lceil6/\tau\rceil$ possible refinements.
The class-count bound gives a uniform bound $K=K(C,\tau,q_0)$
for all partitions in these steps, independently of $n$.
Choose $n_0$ so large that $2CK/\sqrt{n_0}<\tau/2$.
For $n\ge n_0$, all energies in these steps are less than $3$
by \eqref{eq:finite-energy}.

Since the initial energy is nonnegative and each refinement
increases it by more than $\tau/2$, performing all
$\lceil6/\tau\rceil$ steps would give energy greater than $3$,
a contradiction. Thus the process stops with at most $K$ classes.
At that stage, \eqref{eq:averaged-error} gives
$\sum_iw_i\zeta_i<\tau/2+\tau/2=\tau$, as required.
\end{proof}

Apply Claim~2 to the partition from Claim~1.
The resulting partition satisfies \eqref{eq:model-row}, and it
still has $W_5(A)<\tau/2$ because it refines $\mathcal P_0$.
By Claim~1, the initial class count $q_0$ is bounded in terms of
$C,\tau$. Taking the maximum of the bounds from Claim~2 over
these finitely many values of $q_0$ gives $K,n_0$ depending only
on $C,\tau$.
\end{proof}

Lemma~\ref{lem:model}(ii) implies the following bound on
$\|A\|_{\Frob}^2$. Indeed, $A=X+X^{\mathsf T}$ and $X\ge0$ give
$\|A\|_{\Frob}^2\ge2\|X\|_{\Frob}^2$.
Multiplying the row inequalities in \eqref{eq:model-row} by $w_i$
and summing over $i$ yields
\begin{equation}\label{eq:model-energy}
 \tfrac34\|A\|_{\Frob}^2
 \le\|A\|_{\Frob}^2-\tfrac12\|X\|_{\Frob}^2
 \le1+\sum_iw_i\zeta_i<1+\tau.
\end{equation}

\section{Proof of Lemma~\ref{lem:near}}\label{sec:reduction}
We prove Lemma~\ref{lem:near} using the partition from
Lemma~\ref{lem:model}. We first show when lower bounds on the
weighted degrees and the total edge weight force a graph without
triangles or five-cycles to be bipartite. We then apply this
criterion to the graph on the classes and give each original
vertex the color of its class.

\subsection{A bipartiteness criterion for weighted graphs}
For a graph $J$ with positive vertex weights $w_v$, write
$w(S)=\sum_{v\in S}w_v$. The \emph{weighted degree} of $v$ is
$\sigma(v)=w(N(v))$. Write $e_w(J)$ for the \emph{total edge weight},
counting each edge once. Thus
\[
 e_w(J)=\sum_{uv\in E(J)}w_uw_v
       =\tfrac12\sum_{v\in V(J)}w_v\sigma(v).
\]

\begin{lemma}\label{lem:dense}
Let $J$ be a finite simple $\{C_3,C_5\}$-free graph with positive
vertex weights summing to one. If $e_w(J)>2/9$ and
$\sigma(v)\ge1/12$ for every vertex $v$, then $J$ is bipartite.
\end{lemma}
\begin{proof}
We start with an edge whose endpoints have neighborhoods covering
most of the vertex weight. If these neighborhoods cannot be extended
to a bipartition, we will show that the total edge weight is less than $2/9$.
Write $\delta=1/12$. By Cauchy--Schwarz,
\[
 \frac1{e_w(J)}\sum_{uv\in E(J)}w_uw_v\bigl(\sigma(u)+\sigma(v)\bigr)
 =\frac1{e_w(J)}\sum_v w_v\sigma(v)^2
 \ge \frac{(2{e_w(J)})^2}{{e_w(J)}}=4{e_w(J)}>\frac89.
\]
By the above inequality, we can choose an edge $uv$ with $\sigma(u)+\sigma(v)>8/9$. Its endpoints
have no common neighbor, so the set $D=V(J)\setminus(N(u)\cup N(v))$
has weight $r<1/9$.

For each vertex $x$, let $B(x)$ consist of the endpoints of all
two-edge walks starting at $x$. The graph has no closed walk of
length five, since such a walk contains a triangle or a five-cycle.
It follows that $B(x)$ is independent, as an edge joining two such
endpoints would complete a closed walk of length five through $x$.
For the same reason, $L=B(u)$ and $R=B(v)$ are disjoint. Moreover, $N(v)\subseteq L$ and $N(u)\subseteq R$.

Let $Z=V(J)\setminus(L\cup R)$ and $z=w(Z)$. We have $Z\subseteq D$,
and every neighbor of a vertex in $Z$ lies in $D$. Indeed, a neighbor
in $N(u)$ or $N(v)$ would put that vertex of $Z$ in $L$ or $R$.
The set $Z$ is independent. Indeed, the endpoints of an edge in $Z$
would have disjoint neighborhoods by triangle-freeness. Both neighborhoods
would lie in $D$, although their combined weight is at least $2\delta>r=w(D)$.
Thus every edge meeting $Z$ joins $Z$ to $D\setminus Z$.

It remains to rule out a vertex of $Z$ with neighbors on both sides.
Suppose that $x\in Z$ has neighbors in both $L$ and $R$.
Since $N(x)\subseteq D\setminus Z$,
we have $\delta\le w(N(x))\le r-z$, so $z\le r-\delta<\delta$.

We show that the edges between $L$ and $R$ weigh at most
$(1-\delta)^2/4$. A neighbor of $x$ in $R$ has all its neighbors
outside $Z$ in $B(x)\cap L$, and they weigh at least $\delta-z$.
The same holds with the sides reversed. Interchange $L$ and $R$ if
necessary so that
\[
 \delta-z\le w(B(x)\cap L)\le w(B(x)\cap R).
\]
There are no edges between these two sets, since $B(x)$ is independent.

The two sides $L,R$ have combined weight $1-z$. To bring this down
to $1-\delta$, we temporarily remove weight $\delta-z$ from $N(x)$.
Specifically, we multiply the weight of each vertex in $N(x)$ by
$1-(\delta-z)/w(N(x))$. This factor is nonnegative since
$w(N(x))\ge\delta$. The reductions sum to $\delta-z$, and all of
$N(x)$ lies in $L\cup R$. Thus the new side weights sum to
$(1-z)-(\delta-z)=1-\delta$. Their product is the total weight of
all pairs with one vertex on each side. For this fixed sum, the
product is largest when the two weights are equal, so it is at most
$(1-\delta)^2/4$.

Among these pairs, all those between $B(x)\cap L$ and $B(x)\cap R$
are missing. They retain their original weight because triangle-freeness
gives $B(x)\cap N(x)=\varnothing$. Subtracting their weight, we see
that the edges between the sides now weigh at most
\[
 \frac{(1-\delta)^2}{4}-w(B(x)\cap L)\,w(B(x)\cap R).
\]
When we restore weight at a vertex of $N(x)$, the edge weight it adds
between the sides is the restored weight times the weight of its
neighbors on the opposite side. They lie in $B(x)$, so their weights
were not reduced, and their total weight is at most $w(B(x)\cap R)$.
Restoring $\delta-z$ in total
therefore adds at most $(\delta-z)w(B(x)\cap R)$, which is at most
the weight just subtracted since $w(B(x)\cap L)\ge\delta-z$.

Thus the original edge weight between $L$ and $R$ is at most
$(1-\delta)^2/4$. The edges meeting $Z$ weigh at most
$z(r-z)\le r^2/4$, giving
\[
 e_w(J)\le\frac{(1-\delta)^2}{4}+\frac{r^2}{4}
 <\frac14\left(\frac{11}{12}\right)^2+\frac1{324}<\frac29,
\]
contrary to $e_w(J)>2/9$. Thus no vertex of $Z$ has neighbors in both
$L$ and $R$. Since $Z$ is independent, we can place each of its
vertices opposite its neighbors. Hence $J$ is bipartite.
\end{proof}

\subsection{A partition with few spokes inside the parts}
To apply Lemma~\ref{lem:dense}, we need a simple graph on the
classes with no triangles or five-cycles. For a symmetric
nonnegative matrix, its \emph{support} consists of the pairs with
positive entries, including diagonal pairs when present. When the
diagonal is zero, we also regard the support as a simple graph on
the index set.

Lemma~\ref{lem:model} makes $W_5(A)$ small. The following weighted
removal lemma allows us to delete pairs so that the remaining
support has no closed five-step walks. It is the $C_5$ case of
Theorem~4.1 of \cite{CFSZ}, applied with vertex weights $w_i$.
For every $\rho>0$, there is
$\delta>0$ such that any symmetric zero--one matrix $R$ on a finite
probability space with positive weights and $W_5(R)\le\delta$ can be
changed to a matrix $R'$,
by setting entries to zero symmetrically, with
\[
 \sum_{i,j}w_iw_j(R_{ij}-R'_{ij})\le\rho,
 \qquad W_5(R')=0.
\]
The constant $\delta$ depends only on $\rho$, not on the number
or weights of the classes.

Diagonal entries of $R$ may be one. Since all class
weights are positive, $W_5(R')=0$ excludes every closed five-step
walk in its support. This excludes loops and five-cycles. It also
excludes triangles, since a triangle with vertices $a,b,c$ supports
the walk $a,b,a,b,c,a$.
Thus the support of $R'$ is a simple $\{C_3,C_5\}$-free graph on the
same class set.

To apply Lemma~\ref{lem:dense} to the graph on the classes, we use the
minimum triangle degree to bound its weighted degrees from below and
the total triangle count to bound its edge weight from below.

\begin{proof}[Proof of Lemma~\ref{lem:near}]
We first keep only class pairs whose entries in $A$ exceed a small
threshold. These pairs account for almost all spokes. We show that
enough of them remain to apply Lemma~\ref{lem:dense} after we
eliminate closed five-step walks.

Choose $0<\epsilon<\min\{1/4,b,4c\}$ sufficiently small in terms of
$b,c,\eta$. Let $\delta$ be the constant in the weighted removal lemma
for $\rho=\epsilon^2$, and take $0<\tau<\min\{\epsilon^2,\delta\epsilon^5\}$.
By Lemma~\ref{lem:counts}(iv), $\Delta(H)\le C\sqrt n$ for a constant
$C$ depending only on $b$. For all sufficiently large $n$,
Lemma~\ref{lem:model} therefore gives a partition with parameter $\tau$.
By \eqref{eq:model-energy}, $\|A\|_{\Frob}^2\le\tfrac43(1+\tau)<2$.

Let $R_{ij}=1$ if $A_{ij}>\epsilon$ and $R_{ij}=0$ otherwise,
including when $i=j$. By \eqref{eq:model-identities}, the weighted
sum of $X_{ij}+X_{ji}/2$ in each row is at least $b$, and the
weighted sum of all entries of $A$ is at least $4c$. Pairs with
$R_{ij}=0$ contribute at most $\epsilon$ to either sum, since
$X_{ij}+X_{ji}/2\le A_{ij}$. We apply Cauchy--Schwarz to the
retained entries. In the first line, the minimum triangle degree
gives a lower bound on the total weight of the neighbors retained by
class $i$. In the second, the total triangle count gives a lower
bound on the weight of all retained ordered pairs.
\[
\begin{aligned}
 b-\epsilon
 &\le \sum_{j:R_{ij}=1}w_j(X_{ij}+X_{ji}/2)
 \le \left((1+\zeta_i)\sum_{j:R_{ij}=1}w_j\right)^{1/2},\\
 4c-\epsilon
 &\le \sum_{i,j:R_{ij}=1}w_iw_jA_{ij}
 \le \|A\|_{\Frob}\left(\sum_{i,j:R_{ij}=1}w_iw_j\right)^{1/2}.
\end{aligned}
\]
For the first line, $A_{ij}=X_{ij}+X_{ji}$ and $X\ge0$ give
$(X_{ij}+X_{ji}/2)^2\le A_{ij}^2-X_{ji}^2/2$, so
\eqref{eq:model-row} bounds the weighted sum of these squares by
$1+\zeta_i$.

We now delete pairs so that the remaining support has no triangle
or five-cycle. Every term counted by $W_5(R)$ uses five entries of
$A$ larger than $\epsilon$. Thus
$W_5(R)\le\epsilon^{-5}W_5(A)<\delta$ by \eqref{eq:model-walk}.
The weighted removal lemma deletes ordered pairs of total weight
at most $\epsilon^2$ and leaves a simple $\{C_3,C_5\}$-free graph.

Some classes may have lost too many neighbors in this step. Discard
each class $i$ for which the column weights $w_j$ of the entries
just deleted from row $i$ sum to more than $\epsilon$. Averaging
these sums with weights $w_i$ counts each deleted ordered pair
$(i,j)$ with weight $w_iw_j$, giving at most $\epsilon^2$ in total.
Thus these classes have total weight at most $\epsilon$. Also discard
the classes with $\zeta_i>\epsilon$, whose total weight is at most
$\tau/\epsilon<\epsilon$ by \eqref{eq:model-row}.

Let $J$ be the graph on the surviving classes, with their original
weights. In each surviving row the pairs deleted by the removal
lemma cost at most $\epsilon$ in neighbor weight, and discarding
classes costs at most another $2\epsilon$. Ordered pairs with a
discarded first index have total weight at most $2\epsilon$, as do
those with a discarded second index. Counting each edge of $J$ in
both orders, the two Cauchy--Schwarz estimates above therefore give
\[
\begin{aligned}
 w(N_J(i)) \ge \frac{(b-\epsilon)^2}{1+\epsilon}-3\epsilon
                 >\frac1{12}\ \ \text{for}\ i\in V(J), \ \text{and}\ 
 2e_w(J) \ge \frac{(4c-\epsilon)^2}{\tfrac43(1+\tau)}
                  -\epsilon^2-4\epsilon >\frac49.
\end{aligned}
\]
The strict inequalities hold for sufficiently small $\epsilon$ since
$b^2>1/12$, $12c^2>4/9$ and $\tau<\epsilon^2$.
Dividing the surviving weights by their sum can only increase the
weighted degrees and the total edge weight. Lemma~\ref{lem:dense}
now gives a bipartition of $J$.

Give the surviving classes the two colors of $J$, color the discarded
classes arbitrarily, and give each original vertex its class color.
To count spokes, return to the original weights. A class pair $(i,j)$
has weight $w_iw_j$, and $n^{3/2}w_iw_jA_{ij}$ counts the ordered
pairs of vertices in $V_i\times V_j$ joined by a spoke.
If a spoke has both endpoints of the same color and $R_{ij}=1$ for
their classes, then the removal lemma deleted the pair or one of
its classes was discarded. Such pairs have total ordered weight
at most $\epsilon^2+4\epsilon$, so their contribution to
$\sum_{i,j}w_iw_jA_{ij}$ is at most
$\|A\|_{\Frob}\sqrt{\epsilon^2+4\epsilon}$ by Cauchy--Schwarz.
The pairs with $A_{ij}\le\epsilon$ contribute at most $\epsilon$.
Each spoke is counted in both orders, so the number with endpoints
of the same color is at most
\[
 \frac{n^{3/2}}{2}\left(\epsilon+
           \|A\|_{\Frob}\sqrt{\epsilon^2+4\epsilon}\right)
 <\eta n^{3/2}
\]
for sufficiently small $\epsilon$.
\end{proof}

\section*{Acknowledgments}
Guorong Gao received support from the National Key R\&D Program of China
(Grant No.~2023YFA1010202), the National Natural Science Foundation of
China (Grant No.~12401448), and the Natural Science Foundation of Fujian
Province (Grant No.~2024J08030).
Jialin He was supported in part by the National Natural Science
Foundation of China (Grant No.~12601663) and the Science and Technology
Commission of Shanghai Municipality (Grant No.~22DZ2229014).

\section*{Declaration on the use of generative AI}

The authors used ChatGPT-6 Astra to assist in discussing proof strategies, checking proofs, and improving the exposition. 
This project has been under investigation since 2021. Some of the key ideas were contributed by the author.


\begin{thebibliography}{99}
\setlength{\itemsep}{2pt}
\setlength{\parsep}{0pt}
\setlength{\parskip}{0pt}
\bibitem{AS}
N.~Alon and C.~Shikhelman,
\href{https://doi.org/10.1016/j.jctb.2016.03.004}{\emph{Many $T$ copies in $H$-free graphs}},
J. Combin. Theory Ser. B \textbf{121} (2016), 146--172.

\bibitem{BG}
B.~Bollob\'as and E.~Gy\H{o}ri,
\href{https://doi.org/10.1016/j.disc.2007.08.016}{\emph{Pentagons vs. triangles}},
Discrete Math. \textbf{308} (2008), 4332--4336.

\bibitem{CFSZ}
D.~Conlon, J.~Fox, B.~Sudakov and Y.~Zhao,
\href{https://doi.org/10.1112/jlms.12500}{\emph{The regularity method for graphs with few $4$-cycles}},
J. Lond. Math. Soc. \textbf{104} (2021),
2376--2401; arXiv:2004.10180.

\bibitem{Erdos84}
P.~Erd\H{o}s,
\emph{On some problems in graph theory, combinatorial analysis and
combinatorial number theory},
in \emph{Graph Theory and Combinatorics} (Cambridge, 1983),
B.~Bollob\'as (ed.), Academic Press, London, 1984, 1--17.

\bibitem{ES}
P.~Erd\H{o}s and M.~Simonovits,
\href{https://doi.org/10.1007/BF02579234}{\emph{Compactness results in extremal graph theory}},
Combinatorica \textbf{2} (1982), 275--288.

\bibitem{EM}
B.~Ergemlidze and A.~Methuku,
\href{https://doi.org/10.1002/jgt.22723}{\emph{Triangles in $C_5$-free graphs and hypergraphs of girth six}},
J. Graph Theory \textbf{99} (2022), 26--39; arXiv:1811.11873.

\bibitem{EGMS}
B.~Ergemlidze, A.~Methuku, N.~Salia and E.~Gy\H{o}ri,
\href{https://doi.org/10.1002/jgt.22390}{\emph{A note on the maximum number of triangles in a $C_5$-free graph}},
J. Graph Theory \textbf{90} (2019), 227--230; arXiv:1706.02830.

\bibitem{FO}
Z.~F\"uredi and L.~\"Ozkahya,
\href{https://doi.org/10.1016/j.dam.2016.10.013}{\emph{On $3$-uniform hypergraphs without a cycle of a given length}},
Discrete Appl. Math. \textbf{216} (2017), 582--588.

\bibitem{GP}
D.~Gerbner and C.~Palmer,
\href{https://doi.org/10.37236/14563}{\emph{Survey of generalized Tur\'an problems---counting subgraphs}},
Electron. J. Combin., Dynamic Surveys (2026), DS27.

\bibitem{Grzesik}
A.~Grzesik,
\href{https://doi.org/10.1016/j.jctb.2012.04.001}{\emph{On the maximum number of five-cycles in a triangle-free graph}},
J. Combin. Theory Ser. B \textbf{102} (2012), 1061--1066.

\bibitem{GL}
E.~Gy\H{o}ri and H.~Li,
\href{https://doi.org/10.1017/S0963548311000629}{\emph{The maximum number of triangles in $C_{2k+1}$-free graphs}},
Combin. Probab. Comput. \textbf{21} (2012), 187--191.

\bibitem{HHKNR}
H.~Hatami, J.~Hladk\'y, D.~Kr\'a\v{l}, S.~Norine and A.~Razborov,
\href{https://doi.org/10.1016/j.jcta.2012.12.008}{\emph{On the number of pentagons in triangle-free graphs}},
J. Combin. Theory Ser. A \textbf{120} (2013), 722--732.

\bibitem{LHL}
Z.~Lv, Z.~He and M.~Lu,
\href{https://doi.org/10.1016/j.aam.2024.102740}{\emph{Many triangles in $C_5$-free graphs}},
Adv. Appl. Math. \textbf{159} (2024), 102740.

\bibitem{Pascadi}
A.~Pascadi,
\href{https://doi.org/10.37236/13400}{\emph{The sparse regularity method with Schatten norms and entropy}},
Electron. J. Combin. \textbf{33}(1) (2026), P1.26.

\bibitem{RS}
I.~Z.~Ruzsa and E.~Szemer\'edi,
\emph{Triple systems with no six points carrying three triangles},
in \emph{Combinatorics}, Vol.~II (Keszthely, 1976),
Colloquia Mathematica Societatis J\'anos Bolyai \textbf{18},
North-Holland, Amsterdam--New York, 1978, 939--945.
\end{thebibliography}
\end{document}